\documentclass[12pt]{article}
\usepackage[english]{babel}
\usepackage{amsthm,amsmath,amssymb,amsfonts,latexsym}
\newtheorem{theorem}{\qquad Theorem}

\newtheorem{corollary}{\qquad Corollary}
\newtheorem{lemma}{\qquad Lemma}

\newtheorem{remark}{\qquad Remark}

\title{On the Bound of Inverse Images of a Polynomial Map
\footnote{The research was supported by the Russian Foundation for Sciences (project
N 14-11-00433).}
}

\author{Ilya Vyugin}

\date{}

\begin{document}

\maketitle





\begin{abstract}
Let $f_1(x),\ldots,f_n(x)$ be some polynomials. The upper bound on the number of $x\in\mathbb F_p$ such that $f_1(x),\ldots,f_n(x)$ are roots of unit of order $t$ is obtained.
This bound generalize the bound of the paper \cite{V-S} to the case of polynomials of degrees greater than one. The bound is obtained over fields of positive characteristic and over the complex field.
\end{abstract}


\section{Introduction}

Consider the field $\mathbb{F}_p=\mathbb{Z}/p\mathbb{Z}$ where $p$ is a prime number. Let $\mathbb{F}_p^*$ be the multiplicative group of the field $\mathbb{F}_p$ and let $\Gamma$ be a  subgroup of $\mathbb{F}_p^*$ of an order $t=|\Gamma|$. Garcia and Voloch have proved that for any subgroup $\Gamma\subseteq\mathbb{F}_p^*,$ such that $|\Gamma|<(p-1)/((p-1)^{\frac{1}{4}}+1)$ and for any $\mu\in\mathbb{F}_p^*$ the inequality
\begin{eqnarray}\label{GV-est}
|\Gamma\cap (\Gamma+\mu)|\leqslant 4|\Gamma|^{2/3}
\end{eqnarray}
holds. Heath-Brown and Konyagin have re-proved this result (see \cite{H-K}). They also have improved it for the case of a set of equations. Shkredov and Vyugin have generalized the bound to the case of several additive shifts (see \cite{V-S}).

\begin{theorem}[Shkredov and V. \cite{V-S}]\label{th:Vyu_Sh_add_shifts}
Let $\Gamma\subseteq\mathbb{F}_p^*$ be a subroup and let $\mu_1,\ldots,\mu_n \in \mathbb{F}_p^{*}$ be pairwise distinct non-zero elements of $\mathbb{F}_p$, $n\geqslant 2$. Suppose that
$$
32n2^{20n\log (n+1)}\leqslant |\Gamma|,\qquad 4n|\Gamma|(|\Gamma|^{\frac{1}{2n+1}}+1)\leqslant p.
$$
Then we have
$$
|\Gamma\cap(\Gamma+\mu_1)\cap\ldots\cap(\Gamma+\mu_n)|\leqslant 4(n+1)(|\Gamma|^{\frac{1}{2n+1}}+1)^{n+1}.
$$
\end{theorem}
In other words this theorem gives us that
$$
|\Gamma\cap(\Gamma+\mu_1)\cap\ldots\cap(\Gamma+\mu_n)\ll_n |\Gamma|^{\frac{1}{2}+\alpha_n},
$$ 
where $1\ll_n |\Gamma|\ll_n p^{1-\beta_n}$, and $\{\alpha_n\},\{\beta_n\}$ are real sequences, such that $\alpha_n,\beta_n\to 0$, $n\to\infty$.

Consider the map 
\begin{eqnarray}\label{map-linear}
x\stackrel{\varphi}{\rightarrow} (x,x-\mu_1,\ldots,x-\mu_n).
\end{eqnarray}
It is easy to see that 
$$
\Gamma\cap(\Gamma+\mu_1)\cap\ldots\cap(\Gamma+\mu_n)=\varphi(\Gamma^{n+1}),
$$ 
where $\Gamma^{n+1}=\Gamma\times\ldots\times\Gamma$ ($(n+1)$-times). 
We generalize Theorem \ref{th:Vyu_Sh_add_shifts} to the case of polynomial map.

Consider a subgroup $\Gamma\subset\mathbb{F}_p^*$, cosets $\Gamma_1,\dots,\Gamma_n$ by subgroup $\Gamma$ ($\Gamma_i=g_i\Gamma$, $g_i\in \mathbb{F}_p^*$) and a map
\begin{equation}\label{map}
f\,:\, x\longmapsto (f_1(x),\ldots,f_n(x)),\qquad n\geqslant 2
\end{equation}
with polynomials $f_1(x),\dots,f_n(x)\in\mathbb{F}_p[x]$. 

Let us call the set of polynomials  
\begin{eqnarray}\label{set-poly}
f_1(x),\ldots,f_n(x)
\end{eqnarray}
{\it the admissible set of polynomials}, if there exist such $x_1,\ldots,x_n$ that:
$$
f_i(x_i)=0, \quad f_j(x_i)\not= 0, \quad j\not= i, \quad 0\leqslant i,j\leqslant n,
$$
and $f_i(0)\not= 0$, $i=1,\ldots,n$.

Suppose that polynomials $f_1(x),\ldots,f_n(x)$ have degrees $m_i$ ($\deg f_i(x)=m_i$, $i=1,\dots,n$). Let define the set
$$
M=\{ x\mid f_i(x)\in \Gamma_{i},\,\, i=1,\dots,n\}.
$$
Theorem \ref{th-main} gives us the upper bound of cardinality of $M$.

\begin{theorem}\label{th-main}
Let $\Gamma$ be a subgroup of $\mathbb{F}_p^*$ ($p$ is prime), let $\Gamma_{1},\dots,\Gamma_{n}$ be cosets by subgroup $\Gamma$, $n\geqslant 2$ and let $f_1(x),\ldots,f_n(x)$ be an admissible set of poynomials of degrees $m_1,\ldots,m_n$. Let us suppose that
$$
C_1(\textbf{m}, n)<|\Gamma|<C_2(\textbf{m}, n)p^{1-\frac{1}{2n+1}},
$$
where $C_1(\textbf{m},n),C_2(\textbf{m},n)$ are constants depending only on $n$ and $\textbf{m} = (m_1, \ldots , m_n)$.
Then we have the bound
$$  
|M| \leqslant  C_3(\textbf{m}, n)|\Gamma|^{\frac{1}{2} + \frac{1}{2n}},
$$
where $C_3(\textbf{m},n)$ is a constant which depends on $n$ and $\textbf{m}$.
\end{theorem}

\begin{remark}\label{rem1}
Constants $C_1(\textbf{m},n)$, $C_2(\textbf{m},n)$, $C_3(\textbf{m},n)$ can be given as follow:
$$
C_1(\textbf{m}, n)=2^{2n}m_n^{4n},
\qquad
C_2(\textbf{m}, n)=(n+1)^{-\frac{2n}{2n+1}}(m_1\ldots m_n)^{-\frac{2}{2n+1}},
$$
$$
C_3(\textbf{m}, n)=4(n+1)\mathcal{M}_{n}(m_1\dots m_n)^{\frac{1}{n}}.
$$
\end{remark}

\section{Proof of Theorem \ref{th-main}}

We use the idea of proof of papers~\cite{V-S,VSS}. Let us describe Stepanov method (see \cite{H-K}) applied to a polynomial map.
Let us denote by $M'$ the set:
$$
M'=\{x\mid x\in M,\, xf_1(x)\ldots f_n(x)\not=0\}.
$$
Let us construct such polynomial $\Psi(x)$ that satisfy to the following conditions:

\bigskip

1) if $x\in M'$ then $x$ is a root of the polynomial $\Psi(x)$ of an order not less than $D$;

2) $\Psi(x)\not\equiv 0$. 

\bigskip

If such polynomial $\Psi(x)$ exists, than we have the bound:
\begin{eqnarray}\label{est}
 |M| \leqslant 1+\sum_{i=1}^n m_i+\frac{\deg \Psi (x)}{D}.
\end{eqnarray}
Actually, condition 2) gives us that $\Psi(x)$ is a non-zero polynomial. Condition 1) gives us that $|M'|\leqslant \frac{\deg \Psi (x)}{D}$ because any element of $M'$ is a root of polynomial $\Psi(x)$ having the order not less than $D$ by condition 1). A cardinality of the set $M\setminus M'$ is not greater than $1+\sum_{i=1}^n m_i$ by definition of $M'$. 

Let us construct the polynomial $\Psi(x)$. Let $\mathbf{b}$ be a vector $(b_1,\ldots,b_n)$, and $t=|\Gamma|$. Consider the polynomial
\begin{eqnarray}\label{Psi}
\Psi(x) = \sum_{a,\mathbf{b}} \lambda_{a,\mathbf{b}}x^a f_1^{b_1 t}(x)\cdots f_n^{b_n t}(x),\quad a<A,\,\, b_i<B_i,\,\, i=1, \ldots, n
\end{eqnarray}
with coefficients $\lambda_{a,\mathbf{b}}$. Let us define coefficients $\lambda_{a,\mathbf{b}}$ which satisfy to conditions:
\begin{eqnarray}\label{deriv}
\frac{d^k}{dx^k}\Psi(x) \Bigr|_{x \in M'}=0, \quad\quad k=0, \ldots, D-1.
\end{eqnarray}
Let us supose that $x\in M'$, then the condition (\ref{deriv}) is equivalent to
\begin{eqnarray}\label{deriv-mon}
\Bigl[f_1(x) \ldots f_n(x)\Bigr]^k \frac{d^k}{dx^k}\Psi(x)\Bigr|_{x \in M'}=0,\qquad k=0, \ldots, D-1.
\end{eqnarray}
Let us introduce the polynomial $P_{k,a,\mathbf{b}}(x)$ such that
\begin{eqnarray}\label{eq-fPsi}
[f_1(x)\cdots f_n(x)]^k\frac{d^k}{dx^k}\Bigl(x^af_1^{b_1 t}(x)\cdots f_n^{b_n t}(x)\Bigr)=\\ f_1^{b_1 t}(x)\cdots f_n^{b_n t}(x)P_{k,a,\mathbf{b}}(x).\nonumber
\end{eqnarray}
It is easy to see that polynomials $P_{k,a,\mathbf{b}}(x)$ are identity zeros or  
\begin{eqnarray}\label{est-deg}
\deg P_{k,a,\mathbf{b}}(x) \leqslant A+\mathcal{M}_nk - k, 
\end{eqnarray}
where $\mathcal{M}_k=\sum_{i=1}^{k}m_{i}$. The degree of polynomial in left hand side of equality (\ref{eq-fPsi}) is equal to $\mathcal{M}_nk+a+\sum_{i=1}^nb_im_it-k$, consequently, (\ref{est-deg}) holds. \\
\\
Let us substitute $f_i^{t}(x)$ by $g_i^t$ in formulas (\ref{Psi}) and (\ref{deriv-mon}). Actually, $f_i^{t}(x)=g_i^t,\quad i=1,\dots,n$ if $x\in M$, because if $x\in M$ than $f_i(x)\in \Gamma_i=g_i\Gamma$. It means that $f_i^t(x)=g_i^t$ as $t$ --- is the order of subgroup $\Gamma$.
Consequently, we have the equality:
$$
\Bigl[f_1(x) \ldots f_n(x)\Bigr]^k \frac{d^k}{dx^k} \Bigl(x^{a}f_1^{b_1 t}(x)\cdots f_n^{b_n t}(x)\Bigr) \Bigl |_{x \in M}=(g_1^{b_1t}\dots g_n^{b_nt})P_{k,a,\mathbf{b}}(x),
$$
and
\begin{eqnarray}\label{P_k}
\Bigl[f_1(x) \ldots f_n(x)\Bigr]^k \frac{d^k}{dx^k}\Psi(x) \Bigl |_{x \in M}=\sum_{a,\mathbf{b}} \lambda_{a,\mathbf{b}} P_{k,a,\mathbf{b}}(x)=P_k(x).
\end{eqnarray}
Formulas (\ref{est-deg}) and (\ref{P_k}) gives us that degrees of polynomials $P_k(x)$ are bounded as follow:
\begin{eqnarray}\label{est-deg-P}
\deg P_{k}(x) \leqslant A+\mathcal{M}_n k - k,\qquad k=0,1,\ldots,D-1.
\end{eqnarray}
For implicity of the condition (\ref{deriv}), it is suficent to find such $\lambda_{a,\mathbf{b}}$ that they are not vanish simultaneously and such that the following condition  
\begin{equation}\label{zero:Pk}
\forall k= 0, \ldots, D-1\qquad P_k(x)\equiv 0
\end{equation}
holds. Let us obtain coefficients $\lambda_{a,\mathbf{b}}$. Coefficients of polynomials $P_k(x)$ are homogeneous linear combinations of coefficients $\lambda_{a,\mathbf{b}}$, it follows from (\ref{P_k}). The condition (\ref{zero:Pk}) is equivalent to the system of homogeneous linear equations. The system of linear homogeneous equations has a non-zero solution if the number of variables $\lambda_{a, \mathbf{b}}$ is grater than the number of equations. Note that the number of $\lambda_{a, \mathbf{b}}$ is equal to $AB_1 \ldots B_n$, but the number of equations is equal to the number of coefficients fo all polynomials $P_k(x)$, $k=0,1,\ldots, D-1$. It does not exceed of $AD+\mathcal{M}_n\frac{D^2}{2}$, because there exist  bounds (\ref{est-deg-P}).
Consequently, we have the sufficient condition:
\begin{eqnarray}\label{cond-var}
AD+\mathcal{M}_n\frac{D^2}{2} < AB_1 \ldots B_n,
\end{eqnarray}
of existence of non-zero set $\lambda_{a, \mathbf{b}}$ such that the condition 1) holds.

If $\Psi(x)$ is not identity vanish, than
\begin{eqnarray}\label{est-omega}
|M'| \leqslant \frac{\deg \Psi(x)}{D}.
\end{eqnarray}
We prove that $\Psi(x)\not\equiv 0$ if we prove that products 
\begin{eqnarray}\label{product}
x^af_1^{b_1 t}(x)\ldots f_n^{b_n t}(x)
\end{eqnarray}
where $a<A$, $b_i<B_i$, $i=1,\ldots,n$ are linearly independent over $\mathbb{F}_p$, because $\Psi(x)$ is a linear combination of products (\ref{product}) with coefficients $\lambda_{a,\mathbf{b}}$ ($\lambda_{a,\mathbf{b}}$ do not vanish simultaneously).

\begin{lemma}\label{lemma-indep}
Products
\begin{equation}\label{poly}
x^af_1^{b_1 t}(x)\ldots f_n^{b_n t}(x),\qquad a<A,\quad b_i < B_i,\quad i=1,\ldots, n
\end{equation}
are linearly independent over the field $\mathbb{F}_p$, if
\begin{equation}\label{degree}
A-1+(B_1m_1+\ldots+B_{n-1}m_{n-1}-\mathcal{M}_{n-1}) t<p,
\end{equation}
and
\begin{equation}\label{ord-t}
t>AB_1\dots B_{n-1}+\mathcal{M}_n\frac{(B_1\dots B_{n-1})^2}{2}.
\end{equation}
\end{lemma}

{\it Proof.} Let us prove Lemma \ref{lemma-indep} by induction on $n$. In the case $n=0$ Lemma \ref{lemma-indep} is obvious. Actually, the statement of Lemma \ref{lemma-indep} is equivalent to the condition that the monomials $1,x,\ldots,x^{A-1}$
are linearly independent over $\mathbb{F}_p$. Let us prove a step of induction. Suppose that the products
\begin{eqnarray}\label{syst-prod}
x^{a}f_1^{b_1 t}(x)\ldots f_{n-1}^{b_{n-1} t}(x),\qquad a<A,\quad b_i<B_i,\quad i=1, \ldots n-1,
\end{eqnarray}
are linearly independent. We will prove the step of induction from the case $n-1$ to the case $n$ by contradiction. 

Let us suppose that products (\ref{poly}) are linearly dependent. Then there exists a non-trivial set of coefficients $\widetilde{\lambda}_{a,\mathbf{b}}$, such that
\begin{eqnarray*}
\widetilde{\Psi}(x)=\sum \widetilde{\lambda}_{a,\mathbf{b}} x^a f_1^{b_1 t}(x)\ldots f_n^{b_n t}(x) \equiv 0.
\end{eqnarray*}
Without loss of generality let us suppose that $\min_{a,\mathbf{b}}b_n=0$. If $\min_{a,\mathbf{b}}b_n\not=0$ than consider the polynomial $\widetilde{\Psi}(x)/f_n^{t\min_{a,\mathbf{b}}b_n}(x)$ instead of $\widetilde{\Psi}(x)$.
Let us present the polynomial $\widetilde{\Psi}(x)$ in the following form:
\begin{eqnarray}\label{zero-comb2}
\widetilde{\Psi}(x)=f_n^{t}(x)\sum_{a,\mathbf{b}:b_n \ne 0}\widetilde{\lambda}_{a,\mathbf{b}} x^a f_1^{b_1 t}(x)\ldots f_n^{(b_{n}-1) t}(x)+\\ +\sum_{a,\mathbf{b}:b_n = 0}\widetilde{\lambda}_{a,\mathbf{b}} x^a f_1^{b_1 t}(x)\ldots f_{n-1}^{b_{n-1} t}(x)\equiv 0.\nonumber
\end{eqnarray}
Consider the polynomial
$$
\Phi(x)=\sum_{a,\mathbf{b}: b_n = 0}\widetilde{\lambda}_{a,\mathbf{b}} x^{a}f_1^{b_1 t}(x)\ldots f_{n-1}^{b_{n-1} t}(x).
$$ 
Polynomial $\Phi(x)$ is devided by $f_n^{t}(x)$, because the equality (\ref{zero-comb2}) holds. The first term in equality (\ref{zero-comb2}) is divided by $f_n^{t}(x)$, and all sum is divided by $f_n^{t}(x)$ too, consequently, the second term is divided by $f_n^{t}(x)$.
By the proposition of induction $\Phi\not\equiv 0$.
Consequently, to obtain the contradiction to proposition of induction just to prove that if $f_n^t(x)\mid\Phi(x)$ then $\Phi(x)\equiv 0$.

Rewrite $\Phi(x)$ in the following form:
$$
\Phi(x)=\sum_{\mathbf{b}:b_n=0}H_{\mathbf{b}}(x)f_1^{b_1 t}(x)\ldots f_{n-1}^{b_{n-1} t}(x),
$$
where $H_{\mathbf{b}}(x)=\sum_{a}\widetilde{\lambda}_{a,\mathbf{b}} x^{a}$, and all $\mathbf{b}$ are pairwise distinct, $b_i\in \{ 0,\ldots, B_i-1\}$, $i=1, \ldots, n-1$. Note that for all $\mathbf{b}$: $\deg H_{\mathbf{b}}(x)< A$.

Let us introduce the polynomials $Q_{\tilde{\mathbf{b}}}(x)=H_{\tilde{\mathbf{b}}}(x)f_1^{b_1 t}(x)\ldots f_{n-1}^{b_{n-1} t}(x),$ where $\tilde{\mathbf{b}}$ is the vector $(b_1,\ldots,b_{n-1})$. Polynomials $Q_{\tilde{\mathbf{b}}}(x)$ are linearly independent, because by the proposition of induction the polynomials (\ref{syst-prod}) are linearly independent, and $Q_{\tilde{\mathbf{b}}}(x)$ is the linear combination of polynomials (\ref{syst-prod}).
Let us define: $\mathcal{B}_{n}:=\prod_{i=1}^{n}B_{i}.$

Let us consider the Wronskian
$$
W(x)=\begin{vmatrix}
Q_{(0,\ldots,0)}(x) &\ldots & Q_{(B_1-1,\ldots,B_{n-1}-1)}(x)\\
Q_{(0,\ldots,0)}^{'}(x) & \ldots & Q_{(B_1-1,\ldots,B_{n-1}-1)}^{'}(x)\\
\vdots & \ddots & \vdots\\
Q_{(0,\ldots,0)}^{(\mathcal{B}_{n-1}-1)}(x) &\ldots & Q_{(B_1-1,\ldots,B_{n-1}-1)}^{(\mathcal{B}_{n-1}-1)}(x)
\end{vmatrix}.
$$
It is constructed by functions $Q_{\tilde{\mathbf{b}}}(x)$, $\mathbf{b}=(b_1,\ldots,b_n)$, $b_i=0,1,\ldots,B_i$, $i=1,\ldots,n$.  Wronskian $W(x)$ is not equal to zero identity. Actually, the theorem of F.K. Shmidt (see \cite{Shm},\cite{GV-Wron}) states that if Wronskian $W(x)$ is identity vanish then polynomials $Q_{\tilde{\mathbf{b}}}(x)$ are linearly independent over the ring $\mathbb{F}_p[[x^p]]$ of formal power series of the variable $x^p$. The inequality (\ref{degree}) gives us that degrees of polynomials $Q_{\mathbf{b}}(x)$ are less than $p$. It means that the linear independence of $Q_{\tilde{\mathbf{b}}}(x)$ over the ring $\mathbb{F}_p[[x^p]]$ follows from linear independence of $Q_{\tilde{\mathbf{b}}}(x)$ over the field $\mathbb{F}_{p}$.

Wronskian $W(x)$ is devided by
$$
R(x)=\prod_{\mathbf{b}}f_1^{b_1t - \mathcal{B}_{n-1} + 1}(x)\cdot\ldots\cdot f_{n-1}^{b_{n-1}t - \mathcal{B}_{n-1} + 1}(x),
$$ 
because all elements of column with index $\tilde{\textbf{b}}$ are devided by $f_1^{b_1t - \mathcal{B}_{n-1} + 1}(x)\cdot\ldots\cdot f_{n-1}^{b_{n-1}t - \mathcal{B}_{n-1} + 1}(x)$ for each $\textbf{b}$.

Easy to obtain that
$$
\deg (W(x)/R(x)) \leqslant A\mathcal{B}_{n-1} + \mathcal{M}_{n-1}\frac{\mathcal{B}_{n-1}^2}{2}.
$$
By means of elementary tranformations (adding one column to other with some coefficient) Wronskii matrix can be transformed to the form such that elements of one column are function $\Phi(x)$ and its derivatives of orders $1,\ldots,\mathcal{B}_{n-1}-1$.

We have that $(x-x_n)^t\mid \Phi(x)$, because $f_n^t(x)\mid \Phi(x)$, conseqently, we have that $(x-x_n)^{t-(\mathcal{B}_{n-1}-1)}(x)\mid W(x)$. It means that the degree of $(x-x_n)^{t-(\mathcal{B}_{n-1}-1)}(x)$ must be greater than or equal to the degree of $W(x)/R(x)$. It is equivalent to:
$$
t-(\mathcal{B}_{n-1}-1) \leqslant A\mathcal{B}_{n-1} + \mathcal{M}_{n-1}\frac{\mathcal{B}_{n-1}^2}{2}.
$$ 
Consequently, if
\begin{eqnarray}\label{ineq-t}
t\geqslant A\mathcal{B}_{n-1} + \mathcal{M}_{n-1} \frac{\mathcal{B}_{n-1}^2}{2}+\mathcal{B}_{n-1},
\end{eqnarray}
then polynomials (\ref{poly}) are linearly independent. We have proved the step of induction. $\Box$


\subsection{Setting of Parameters}

To prove Theorem \ref{th-main} we have to set the parameters $A,B_1,\ldots,B_n,D$, and proved that they are satisfy to the necessary conditions (\ref{cond-var}), (\ref{degree}), (\ref{ord-t}). The bound can be obtained by substituting of parameters to formulas (\ref{est}), (\ref{Psi}).

Without loss of generality, let us set the following
\begin{eqnarray}\label{m1-n}
m_1\leqslant\ldots\leqslant m_n,
\end{eqnarray}
where $m_i=\deg f_i(x)$, $i=1,\ldots,n$. A permutation of polynomials do not change conditions of Theorem~\ref{th-main}.

Let us put $B=(m_1\dots m_n)^{\frac{1}{n}}t^{\frac{1}{2n}}$,
$$
B_i=\left[\frac{B}{m_i}\right],\quad i=1,\dots,n,\qquad A=B_1\dots B_n,\quad D=\left[\frac{A}{\mathcal{M}_n}\right],
$$
where $[\cdot]$ is the integer part of the number.
Let us check conditions (\ref{cond-var}), (\ref{degree}), (\ref{ord-t}).

The condition (\ref{cond-var}) has the form:
$$
AD+\mathcal{M}_n\frac{D^2}{2}\leqslant\frac{1}{\mathcal{M}_n}A^2+\frac{1}{2\mathcal{M}_n}A^2<A^2=AB_1 \ldots B_n.
$$
It is true, because $\mathcal{M}_n\geqslant n\geqslant 2$.

The condition (\ref{degree}) has the form:
$$
\deg\Psi(x)\leqslant A+\sum_{i=1}^n(B_i-1)m_it<A+nBt\leqslant 
$$ 
\begin{eqnarray}\label{ineq16}
\leqslant t^{1/2}+ n(m_1\dots m_n)^{\frac{1}{n}}t^{1+\frac{1}{2n}}\leqslant  (n+1)(m_1\dots m_n)^{\frac{1}{n}}t^{\frac{2n+1}{2n}}<p.
\end{eqnarray}
Actually, $\left[\frac{B}{m_i}\right]\leqslant\frac{B}{m_i}$ and, consequently,
$$
A=\prod_{i=1}^n\left[\frac{B}{m_i}\right]\leqslant \frac{B^n}{m_1\dots m_n}=t^{1/2}.
$$
The last inequality in (\ref{ineq16}) follows from  
\begin{eqnarray}\label{ineq-t-h}
t<(n+1)^{-\frac{2n}{2n+1}}(m_1\ldots m_n)^{-\frac{2}{2n+1}}p^{1-\frac{1}{2n+1}}.
\end{eqnarray}
The inequality (\ref{degree}) is proved.

Let us show that the condition (\ref{ord-t}) is also holds. Let us consider the right hand side of the inequality (\ref{ord-t}):
$$
AB_1\ldots B_{n-1}+\frac{\mathcal{M}_n}{2}(B_1\ldots B_{n-1})^2=(B_1\ldots B_{n-1})^2\left(B_n+\frac{\mathcal{M}_n}{2}\right) <
$$
(Use now that $B_i=\left[\frac{B}{m_i}\right]\geqslant \gamma\frac{B}{m_i}$, $i=1,\ldots,n$, where $\gamma$ is the following
$$
\gamma=\frac{B-\max_{1\leqslant i\leqslant n}m_i}{B}=\frac{B-m_n}{B},
$$ 
using the condition (\ref{m1-n}))
$$
<\gamma^{2n-1}\frac{(m_1\dots m_n)^{2-2/n}}{(m_1\dots m_{n-1})^{2}}t^{1-1/n}\left(\frac{(m_1\ldots m_n)^{1/n}}{m_n}t^{\frac{1}{2n}}+\frac{\mathcal{M}_{n}}{2}\right).
$$
(Use that $\frac{(m_1\ldots m_n)^{1/n}}{m_n}\leqslant 1$, and that $t>\left(\mathcal{M}_n/2\right)^{2n}$. )
$$
<2\gamma^{2n-1}m_n^2t^{1-\frac{1}{2n}}<t.
$$
The last inequality follows from
$$
t^{\frac{1}{2n}}>2\gamma^{2n-1}m_n^2=2\left(\frac{B-m_n}{B}\right)^{2n-1}m_n^2,
$$
this follows from $\left(\frac{B-m_n}{B}\right)<1$, and 
\begin{eqnarray}\label{ineq-t}
t>2^{2n}m_n^{4n}.
\end{eqnarray}
Use (\ref{est}) for estimation of $|M|$:
$$
|M|\leqslant 1+\mathcal{M}_{n}+\frac{A+\sum_{i=1}^n (B_i-1)m_it}{D}<
$$
(Use now that $t^{1/2}>1+\mathcal{M}_{n}$. It follows from (\ref{ineq-t}).)
$$
<\frac{A+nBt}{D}=\frac{A+nBt}{[A/\mathcal{M}_{n}]}<
$$
(Obtain the inequatity by means (\ref{ineq16}).)
$$
\leqslant \frac{(n+1)(m_1\dots m_n)^{\frac{1}{n}}t^{\frac{2n+1}{2n}}}{\gamma^n t^{1/2}/\mathcal{M}_{n}}\leqslant
$$
(It is easy to see that $\gamma^n=\left(1-\frac{m_n}{B}\right)^n>1/4$, because $\frac{m_n}{B}<1/n$ follows from (\ref{ineq-t}).)
$$
\leqslant 4(n+1)\mathcal{M}_{n}(m_1\dots m_n)^{\frac{1}{n}}t^{\frac{1}{2}+\frac{1}{2n}}.
$$
It is easy to see that constants $C_1(\textbf{m}, n)$, $C_2(\textbf{m}, n)$ и $C_3(\textbf{m}, n)$ can be setted as follows:
from the inequality (\ref{ineq-t}) let us obtain:
$$
C_1(\textbf{m}, n)=2^{2n}m_n^{4n},
$$
let us remind that $m_n=\max_{1\leqslant i\leqslant n} m_i$;
from the inequality (\ref{ineq-t-h}) obtain:
$$
C_2(\textbf{m}, n)=(n+1)^{-\frac{2n}{2n+1}}(m_1\ldots m_n)^{-\frac{2}{2n+1}};
$$
and from the final bound obtain the value of the last constant:
$$
C_3(\textbf{m}, n)=4(n+1)\mathcal{M}_{n}(m_1\dots m_n)^{\frac{1}{n}}.
$$
These constants prove Remark~\ref{rem1}.
$\Box$

Let us consider the linear map (\ref{map-linear}). Let us obtain Corollary~\ref{Cor1} which is the reslt of the paper \cite{V-S}.

\begin{corollary}\label{Cor1}
Let $\Gamma$ be a subgroup of $\mathbb{F}_p^*$ ($p$ is a prime number), $n\geqslant 1$. Let the following inequality 
$$
C_1(n)<|\Gamma|<C_2(n)p^{1-\frac{1}{2n+1}},
$$
where $C_1(n),C_2(n)$ are constant depending only on $n$, holds. Then we have the following bound:
$$  
|\Gamma\cap(\Gamma+\mu_1)\cap\ldots\cap(\Gamma+\mu_{n-1})| \leqslant  C_3(n)|\Gamma|^{\frac{1}{2} + \frac{1}{2n}},
$$
where $C_3(n)$ is some constant depending only on от $n$, holds.
\end{corollary}

\section{Polynomial Maps over $\mathbb{C}$}

Let us consider the analog Theorem~\ref{th-main} for the complex field. Let $G=\{ x\mid x^t=1\}$ be a subgroup of roots of orders $t$ of unity of the group $\mathbb{C}^{*}$. Let us denote cosets of the subgroup $G$ by $G_1,\dots,G_n$. Consider the map 
\begin{equation}\label{map}
f\,:\, x\longmapsto (f_1(x),\ldots,f_n(x)),\qquad n\geqslant 2,
\end{equation}
where $f_1(x),\dots,f_n(x)\in\mathbb{C}[x]$ are polynomials. The definition of admissibility of polynomials is analogous to the definition for polynomials over $\mathbb{F}_p$. 

For the cardinality of the set
$$
M=\{ x\mid f_i(x)\in G_{i},\,\, i=1,\dots,n\}
$$
the following theorem holds.

\begin{theorem}\label{th-comp}
Let $G$ be a subgroup of $\mathbb{C}^*$ of roots of unity of some order, $G_1,\dots,G_n$ are cosets of $G$, $n\geqslant 2$, $f_1(x),\ldots,f_n(x)$ is an admissible set of polynomials of degrees $m_1,\ldots,m_n$. Let us suppose that:
$$
|G|>\tilde{C}_1(\textbf{m}, n),
$$
where $\tilde{C}_1(\textbf{m},n))$ is a constant depending only on $n$ and $\textbf{m}$.
Then we have the following bound:
$$  
|M| \leqslant  \tilde{C}_2(\textbf{m}, n)|\Gamma|^{\frac{1}{2} + \frac{1}{2n}},
$$
where $\tilde{C}_2(\textbf{m},n)$ is a constant depending only on $n$ and $\textbf{m}$.
\end{theorem}

\begin{remark}\label{rem2}
Constants can be setted as follows: $\tilde{C}_1(\textbf{m},n)=C_1(\textbf{m},n)$, $\tilde{C}_2(\textbf{m},n)=C_3(\textbf{m},n)$.
\end{remark}

The proofs of Theorem \label{th-comp} and Remark \ref{rem2} almost completely repeat the proofs of Theorem~\ref{th-main} and Remark~\ref{rem1}. We will not repeat these proofs. We only describe two small changes. We do not reqire that degree of polynomial $\Psi(x)$ is the less than characteristic of the field. It gives us that the restriction (\ref{ineq16}) is not actual. Also instead of theorem of F.K. Shmidt we use the theorem on linear dependence of a set of functions and vanishing of Wronskian.

Vyugin I.V.\\
Insitute for Information Transmission Problems RAS,  \\
and\\
National Research University Higher School of Economics,\\
{\it vyugin@gmail.com}.\\

\end{document}